\documentclass[11pt,reqno]{amsart}
\usepackage{mathrsfs}
\usepackage{amssymb}
\usepackage{amsmath}
\usepackage{amscd}
\usepackage[colorlinks, linkcolor=blue, anchorcolor=blue, citecolor=blue]{hyperref}
\begin{document}
\setlength{\oddsidemargin}{0cm} \setlength{\evensidemargin}{0cm}
\baselineskip=20pt

\theoremstyle{plain} \makeatletter
\newtheorem{theorem}{Theorem}[section]
\newtheorem{proposition}[theorem]{Proposition}
\newtheorem{lemma}[theorem]{Lemma}
\newtheorem{coro}[theorem]{Corollary}

\theoremstyle{definition}
\newtheorem{defi}[theorem]{Definition}
\newtheorem{notation}[theorem]{Notation}
\newtheorem{exam}[theorem]{Example}
\newtheorem{prop}[theorem]{Proposition}
\newtheorem{conj}[theorem]{Conjecture}
\newtheorem{prob}[theorem]{Problem}
\newtheorem{remark}[theorem]{Remark}
\newtheorem{claim}{Claim}

\newcommand{\SO}{{\mathrm S}{\mathrm O}}
\newcommand{\SU}{{\mathrm S}{\mathrm U}}
\newcommand{\Sp}{{\mathrm S}{\mathrm p}}
\newcommand{\G}{{\mathrm G}}
\newcommand{\so}{{\mathfrak s}{\mathfrak o}}
\newcommand{\Ad}{{\mathrm A}{\mathrm d}}
\newcommand{\m}{{\mathfrak m}}
\newcommand{\g}{{\mathfrak g}}

\numberwithin{equation}{section}

\title{Non-naturally reductive Einstein metrics on $\mathrm{SO}(n)$}
\author{Huibin Chen}
\address{School of Mathematical Sciences and LPMC, Nankai University,
Tianjin 300071, P.R. China}\email{chenhuibin@mail.nankai.edu.cn}
\author{Zhiqi Chen}
\address{School of Mathematical Sciences and LPMC, Nankai University,
Tianjin 300071, P.R. China}\email{chenzhiqi@nankai.edu.cn}
\author{Shaoqiang Deng}
\address{School of Mathematical Sciences and LPMC, Nankai University,
Tianjin 300071, P.R. China}\email{dengsq@nankai.edu.cn}
\date{}

\begin{abstract}
In this article, we prove that every compact simple Lie group $\SO(n)$ for $n\geq 10$ admits at least $2\left([\frac{n-1}{3}]-2\right)$ non-naturally reductive left-invariant Einstein metrics.
\end{abstract}

\maketitle

\section{Introduction}
A Riemannian manifold $(M,\langle\cdot,\cdot\rangle)$ is called Einstein if there exists a constant $\lambda$ such that $\mathrm{Ric}=\lambda \langle\cdot,\cdot\rangle$, where $\mathrm{Ric}$ is the Ricci tensor of the metric $\langle\cdot,\cdot\rangle$ . For a survey of the results on Einstein metrics, we refer to the book \cite{Be} and the papers \cite{W1,W2}. For the study on Einstein metrics of homogeneous manifolds, see \cite{Bo,BoWaZi,WaZi}.

There are a lot of important results on Einstein metrics on Lie groups which are a special class of homogeneous manifolds. It is shown in \cite{La} that any Einstein solvmanifold, i.e. a simply connected solvable Lie group admitting a left-invariant metric, is standard. According to a conjecture of Alekseevskii (\cite{Be}), this exhausts noncompact Lie groups with left-invariant Einstein metrics. Also Einstein metrics on solvable Lie groups are unique to isometry and scaling \cite{He}, which is in sharp contrast to the compact setting.

In \cite{JW},  D' Atri and Ziller prove every compact simple Lie group except $\SO(3)$ admits at least two left-invariant Einstein metrics which are all naturally reductive. Meanwhile, they pose the question whether there exist non-naturally reductive left-invariant Einstein metrics on compact simple Lie groups. From then on, there are a lot of studies on non-naturally reductive Einstein metrics on compact simple Lie groups. In \cite{Mo}, Mori obtains non-naturally reductive left-invariant Einstein metrics on $\SU(n)$ for $n\geq6$. Then in \cite{ArMoSa}, the authors prove the existence of left-invariant Einstein metrics on compact Lie groups $\SO(n)$ for $n\geq 11$, $\Sp(n)$ for $n\geq3 $, $\mathrm{E}_6, \mathrm{E}_7$ and $\mathrm{E}_8$. After that, Chen and Liang \cite{ChLi} give one non-naturally reductive left-invariant Einstein metric on  $\mathrm{F}_4$. In \cite{ArSaSt3}, Arvanitoyeorgos, Sakane and Statha obtain left-invariant Einstein metrics on compact Lie groups $\SO(n)$ for $n\geq7$, which are not naturally reductive. Recently, Chrysikos and Sakane find new non-naturally reductive Einstein metrics on exceptional Lie groups in \cite{ChSa}, especially they give the first non-naturally reductive Einstein metric on $\G_2$.

That is to say, every compact simple Lie group except some small groups admits non-naturally reductive left-invariant Einstein metrics. But how many non-naturally reductive left-invariant Einstein metrics does every compact simple Lie group possibly admit? It is natural to discuss the lower bound of the number of non-naturally reductive left-invariant Einstein metrics on compact simple Lie groups, especially classical Lie groups.

Yan and Deng prove an interesting result in \cite{YD1}: for any integer $n=p_1^{l_1}p_2^{l_2}\cdots p_s^{l_s}$ with $p_i$ prime and $p_i\not=p_j$,
\begin{enumerate}
  \item $\SO(2n)$ admits at least $(l_1+1)(l_2+1)\cdots (l_s+1)-3$ non-equivalent non-naturally reductive Einstein metrics,
  \item and $\Sp(2n)$ admits at least $(l_1+1)(l_2+1)\cdots (l_s+1)-1$ non-equivalent non-naturally reductive Einstein metrics.
\end{enumerate}
That is, they give lower bounds of the number for $\SO(2n)$ and $\Sp(2n)$, which depend on the integer $n$. But we point out that the above results hold only when $n$ is not a prime. A much better estimate for $\Sp(n)$ is given in \cite{CC} that every $\Sp(n)$ admits at least $2[\frac{n-1}{3}]$ non-naturally reductive left-invariant Einstein metrics for $n\geq 4$.

In this article, we obtain the following lower bound of the number of non-naturally reductive left-invariant Einstein metrics on $\SO(n)$.

\begin{theorem}\label{thm}
 For every integer $n\geq10$, $\SO(n)$ admits at least $2\left([\frac{n-1}{3}]-2\right)$ non-naturally reductive left-invariant Einstein metrics.
\end{theorem}

The paper is organized as follows. In section 2, we recall the study on non-naturally reductive left-invariant metrics on $\SO(n)$ in \cite{ArSaSt3}, in particular the Ricci tensor of a class of left-invariant metrics $\langle\cdot,\cdot\rangle$ on $\SO(n)$, and sufficient and necessary conditions for $\langle\cdot,\cdot\rangle$ to be naturally reductive. Furthermore, they prove that $\SO(n)$ admits non-naturally reductive Einstein metrics which are $\Ad(\SO(n-6)\times\SO(3)\times\SO(3))$-invariant. In section 3, based on the Ricci tensor formulae in \cite{ArSaSt3} and the technique of Gr\"{o}bner basis, we prove that $\SO(2k+l)$ admits at least two non-naturally reductive Einstein metrics which are $\Ad(\SO(k)\times\SO(k)\times\SO(l))$-invariant when $l>k\geq 3$. It implies Theorem~\ref{thm}.

\section{The study on $\SO(n)$ in \cite{ArSaSt3}}
Let $G$ be a compact semisimple Lie group and let $K$ be a connected closed subgroup of $G$ with Lie algebras $\mathfrak{g}$ and $\mathfrak{k}$. Since the Killing form $B$ of $\mathfrak{g}$ is negative definite, $-B$ is an $\Ad(G)$-invariant inner product on $\mathfrak{g}$. Let $\mathfrak{g}=\mathfrak{k}\oplus\mathfrak{m}$ be the reductive decomposition of $\mathfrak{g}$ such that $[\mathfrak{k},\mathfrak{m}]\subset\mathfrak{m}$. Then $\mathfrak{m}$ can be identified with the tangent space of $G/K$ at the origin. Assume that $\mathfrak{m}$ admits a decomposition into mutually non-equivalent irreducible $\Ad(K)$-modules:
\begin{equation}
\mathfrak{m}=\mathfrak{m}_1\oplus\cdots\oplus\mathfrak{m}_q
\end{equation}
Then any $G$-invariant metric on $G/K$ has the following form
\begin{equation}\label{met3}
\langle\cdot,\cdot\rangle=x_1(-B)|_{\mathfrak{m}_1}+\cdots+x_q(-B)|_{\mathfrak{m}_q},
\end{equation}
where $x_1,\cdots,x_q\in\mathbb{R}^+$. The $G$-invariant metric $\langle\cdot,\cdot\rangle$ on $G/K$ is called naturally reductive if
\begin{equation*}
\langle [X,Y]_{\mathfrak m},Z\rangle+\langle Y,[X,Z]_{\mathfrak m}\rangle=0, \quad \forall X,Y,Z\in\mathfrak m.
\end{equation*}
Note that $G$-invariant symmetric covariant $2$-tensors on $G/K$ are of the same form as Riemannian metrics. In particular, the Ricci tensor $r$ of a $G$-invariant Riemannian metric on $G/K$ is of the same form as (\ref{met3}), that is,
\begin{equation}
r=y_1(-B)|_{\mathfrak{m}_1}+\cdots+y_q(-B)|_{\mathfrak{m}_q},
\end{equation}
where $y_1,\cdots, y_q\in \mathbb{R}$. Let $d_i=\dim\mathfrak m_i$ and let $\{e^{i}_{\alpha}\}_{\alpha=1}^{d_i}$ be a $(-B)$-orthonormal basis of $\mathfrak m_i$. Denote  $A_{\alpha,\beta}^{\gamma}=-B([e_{\alpha}^{i},e_{\beta}^{j}],e_{\gamma}^{k})$, i.e.  $[e_{\alpha}^{i},e_{\beta}^{j}]=\sum_{\gamma}A_{\alpha,\beta}^{\gamma}e_{\gamma}^{k}$, and define
\begin{equation*}
(ijk):=\sum(A_{\alpha,\beta}^{\gamma})^2,
\end{equation*}
where the sum is taken over all indices $\alpha,\beta,\gamma$ with $e_{\alpha}^{i}\in\mathfrak m_i, e_{\beta}^{j}\in\mathfrak m_j, e_{\gamma}^{k}\in\mathfrak m_k$. Then $(ijk)$ is independent of the choice for the $-B$-orthonormal basis of $\mathfrak m_i,\mathfrak m_j,\mathfrak m_k$, and $(ijk)=(jik)=(jki)$. Then we have the following result.

\begin{lemma}[\cite{PaSa}]\label{lem1}
The components $r_1,\cdots,r_{p+q}$ of the Ricci tensor $r$ associated to $\langle\cdot,\cdot\rangle$ of the form~\eqref{met3} on $G/K$
\begin{equation*}
r_k=\frac{1}{2x_k}+\frac{1}{4d_k}\sum_{j,i}\frac{x_k}{x_jx_i}(kji)-\frac{1}{2d_k}\sum_{j,i}\frac{x_j}{x_kx_i}(jki),\quad(k=0,1,\cdots,p+q).
\end{equation*}
Here, the sums are taken over all $i=1,\cdots,p+q$.
\end{lemma}

For $G=\SO(k_1+k_2+k_3)$, we consider the closed subgroup $K=\SO(k_1)\times \SO(k_2) \times \SO(k_3)$, where the embedding of $K$ in $G$ is diagonal, and we have the fibration
\begin{equation*}
\SO(k_1)\times \SO(k_2)\times \SO(k_3)\rightarrow \SO(k_1+k_2+k_3)\rightarrow \SO(k_1+k_2+k_3)/(\SO(k_1)\times \SO(k_2)\times \SO(k_3)).
\end{equation*}
Denote the tangent space of $\SO(k_1+k_2+k_3)$, $\SO(k_1+k_2+k_3)/(\SO(k_1)\times \SO(k_2)\times \SO(k_3))$ and $\SO(k_1)\times \SO(k_2)\times \SO(k_3)$ at the origin by $\so(k_1+k_2+k_3)$, $\m$ and $\so(k_1)\oplus \so(k_2) \oplus \so(k_3)$, respectively. Then we have
\begin{equation*}
\so(k_1+k_2+k_3)=\so(k_1)\oplus \so(k_2)\oplus \so(k_3)\oplus\m.
\end{equation*}
By setting $\m_1=\so(k_1),\m_2=\so(k_2)$ and $\m_3=\so(k_3)$, we have the following decomposition:
\begin{equation*}
\so(k_1+k_2+k_3)=\m_1\oplus \m_2\oplus \m_3\oplus \m_{12}\oplus \m_{13}\oplus \m_{23},
\end{equation*}
where $\m_{12},\m_{13}$ and $\m_{23}$ are irreducible submodules of $\m$.
Note that there is a diffeomorphism
$$G/\{e\}\cong(G\times \SO(k_1)\times \SO(k_2)\times \SO(k_3))/\mathrm{diag}(\SO(k_1)\times \SO(k_2)\times \SO(k_3)).$$
Consider left-invariant metrics on $G$ which are determined by $\Ad(K)$-invariant inner products on $\so(k_1+k_2+k_3)$ given by
\begin{equation}\label{met1}
\begin{split}
\langle\cdot,\cdot\rangle & = x_1(-B)|_{\m_1}+x_2(-B)|_{\m_2}+x_3(-B)|_{\m_3}\\
&+x_{12}(-B)|_{\m_{12}}+x_{13}(-B)|_{\m_{13}}+x_{23}(-B)|_{\m_{23}}.
\end{split}
\end{equation}
By Lemma~\ref{lem1}, we have the following formulae of the Ricci tensor.
\begin{lemma}[\cite{ArSaSt3}]\label{Prop:0}
The components of the Ricci tensor $r$ of left-invariant metrics on $G$ defined by $\langle\cdot,\cdot\rangle$ of the form $(\ref{met1})$ are given as follows:
\begin{equation}
\begin{split}
r_1&=\frac{k_1-2}{4(n-2)x_1}+\frac{1}{4(n-2)}\left(k_2\frac{x_1}{x_{12}^2}+k_3\frac{x_1}{x_{13}^2}\right),\\
r_2&=\frac{k_2-2}{4(n-2)x_2}+\frac{1}{4(n-2)}\left(k_1\frac{x_2}{x_{12}^2}+k_3\frac{x_2}{x_{23}^2}\right),\\
r_3&=\frac{k_3-2}{4(n-2)x_3}+\frac{1}{4(n-2)}\left(k_1\frac{x_3}{x_{13}^2}+k_2\frac{x_3}{x_{23}^2}\right),\\
r_{12}&=\frac{1}{2x_{12}}+\frac{k_3}{4(n-2)}\left(\frac{x_{12}}{x_{13}x_{23}}-\frac{x_{13}}{x_{12}x_{23}}-\frac{x_{23}}{x_{12}x_{13}}\right)-\frac{1}{4(n-2)}\left(\frac{(k_1-1)x_1}{x_{12}^2}+\frac{(k_2-1)x_2}{x_{12}^2}\right),\\
r_{13}&=\frac{1}{2x_{13}}+\frac{k_2}{4(n-2)}\left(\frac{x_{13}}{x_{12}x_{23}}-\frac{x_{12}}{x_{13}x_{23}}-\frac{x_{23}}{x_{12}x_{13}}\right)-\frac{1}{4(n-2)}\left(\frac{(k_1-1)x_1}{x_{13}^2}+\frac{(k_3-1)x_3}{x_{13}^2}\right),\\
r_{23}&=\frac{1}{2x_{23}}+\frac{k_1}{4(n-2)}\left(\frac{x_{23}}{x_{13}x_{12}}-\frac{x_{13}}{x_{12}x_{23}}-\frac{x_{12}}{x_{23}x_{13}}\right)-\frac{1}{4(n-2)}\left(\frac{(k_2-1)x_2}{x_{23}^2}+\frac{(k_3-1)x_3}{x_{23}^2}\right),
\end{split}
\end{equation}
where $n=k_1+k_2+k_3$.
\end{lemma}

\begin{defi}
A Riemannian homogeneous space $(M=G/K, \langle\cdot,\cdot\rangle)$ with a reductive complement $\m$ of $\mathfrak{k}$ in $\mathfrak{g}$ is called naturally reductive if
$$\langle [X,Y]_\m,Z\rangle + \langle Y,[X,Z]_\m\rangle=0,\,\forall X,Y,Z\in\m.$$
\end{defi}
Based on the criterion for a left-invariant metric on a Lie group to be naturally reductive given in \cite{JW}, we have the following lemma.
\begin{lemma}[\cite{ArSaSt3}]\label{prop}
If a left-invariant metric $\langle\cdot,\cdot\rangle$ of the form $(\ref{met1})$ on $G$ is naturally reductive with respect to $G\times L$, where $L$ is a closed subgroup of $G=\SO(k_1+k_2+k_3)$, then one of the following holds:
\begin{enumerate}
\item $x_1=x_2=x_{12}$, $x_{13}=x_{23};$
 \item $x_2=x_3=x_{23}$, $x_{12}=x_{13};$
 \item $x_1=x_3=x_{13}$, $x_{12}=x_{23};$
 \item $x_{12}=x_{13}=x_{23}$.
\end{enumerate}
Conversely, if one of the above conditions is satisfied, then there exists a closed subgroup $L$ of $G$ such that the metric $\langle\cdot,\cdot\rangle$ of the form $(\ref{met1})$ is naturally reductive with respect to $G\times L$.
\end{lemma}

\section{The proof of Theorem~\ref{thm}}
In \cite{ArSaSt3}, the authors prove that, for any $n \geq 9$, the Lie group $\SO(n)$ admits at least one left-invariant Einstein metric determined by the $\Ad(\SO(n-6)\times\SO(3)\times\SO(3))$-invariant inner product of the form~\eqref{met1}, which is non-naturally reductive. In this section, we first prove the following theorem.
\begin{theorem}\label{thm2}
For any $l>k\geq 3$, $\SO(2k+l)$ admits at least two left-invariant Einstein metrics determined by $\Ad(\SO(k)\times\SO(k)\times\SO(l))$-invariant inner products of the form~\eqref{met1}, which are non-naturally reductive.
\end{theorem}

We will prove Theorem~\ref{thm2} by solving homogeneous Einstein equations
$$r_1=r_2,\quad r_2=r_3,\quad r_3=r_{12},\quad r_{12}=r_{13},\quad r_{13}=r_{23}$$
under the assumption $k=k_1=k_2$ and $l=k_3$. Furthermore, we consider the metric \eqref{met1} with $x_1=x_2$. Then we have $x_{13}=x_{23}$. Standard the metric with $x_{13}=x_{23}=1$. Therefore homogeneous Einstein equations are equivalent to the following system of equations:
\[\left\{\begin{aligned}
f_1&=-2\,kx_{{2}}{x_{{3}}}^{2}{x_{{12}}}^{2}+l{x_{{2}}}^{2}x_{{3}}{x_{{12}}
}^{2}+k{x_{{2}}}^{2}x_{{3}}+kx_{{3}}{x_{{12}}}^{2}-lx_{{2}}{x_{{12}}}^
{2}+2\,x_{{2}}{x_{{12}}}^{2}-2\,x_{{3}}{x_{{12}}}^{2}=0,\\
f_2&=l{x_{{2}}}^{2}{x_{{12}}}^{2}-lx_{{2}}{x_{{12}}}^{3}+3\,k{x_{{2}}}^{2}-
4\,x_{{2}}x_{{12}}k+{x_{{12}}}^{2}k-2\,{x_{{2}}}^{2}+4\,x_{{2}}x_{{12}
}-2\,{x_{{12}}}^{2}=0,\\
f_3&=kx_{{2}}x_{{3}}+2\,k{x_{{3}}}^{2}+kx_{{3}}x_{{12}}+l{x_{{3}}}^{2}-4\,k
x_{{3}}-2\,x_{{3}}l-x_{{2}}x_{{3}}-{x_{{3}}}^{2}+l+4\,x_{{3}}-2=0.
\end{aligned}\right.\]
In particular,
$$f_2= \left( x_{{2}}-x_{{12}} \right)  \left( lx_{{2}}{x_{{12}}}^{2}+3\,kx_
{{2}}-kx_{{12}}-2\,x_{{2}}+2\,x_{{12}} \right).$$
If $x_2=x_{12}$, by Proposition \ref{prop}, the left-invariant metric is naturally reductive. Assume that $x_2\neq x_{12}$. Then we have
$$x_2={\frac {x_{{12}} \left( k-2 \right) }{l{x_{{12}}}^{2}+3\,k-2}}.$$
Substitute it into equations $f_1=0$ and $f_3=0$, we have the following equations:
\[\left\{\begin{aligned}
g_1=&2\,{k}^{2}l{x_{{3}}}^{2}{x_{{12}}}^{3}-k{l}^{2}x_{{3}}{x_{{12}}}^{4}-4
\,kl{x_{{3}}}^{2}{x_{{12}}}^{3}+2\,{l}^{2}x_{{3}}{x_{{12}}}^{4}+6\,{k}
^{3}{x_{{3}}}^{2}x_{{12}}-7\,{k}^{2}lx_{{3}}{x_{{12}}}^{2}+k{l}^{2}{x_
{{12}}}^{3}\\&-16\,{k}^{2}{x_{{3}}}^{2}x_{{12}}+20\,klx_{{3}}{x_{{12}}}^{
2}-2\,kl{x_{{12}}}^{3}-2\,{l}^{2}{x_{{12}}}^{3}-10\,{k}^{3}x_{{3}}+3\,
{k}^{2}lx_{{12}}+8\,k{x_{{3}}}^{2}x_{{12}}\\&-12\,lx_{{3}}{x_{{12}}}^{2}+
4\,l{x_{{12}}}^{3}+34\,{k}^{2}x_{{3}}-6\,{k}^{2}x_{{12}}-8\,klx_{{12}}
-32\,kx_{{3}}+16\,kx_{{12}}+4\,x_{{12}}l+8\,x_{{3}}\\&-8\,x_{{12}}=0,\\
g_2=&2\,kl{x_{{3}}}^{2}{x_{{12}}}^{2}+klx_{{3}}{x_{{12}}}^{3}+{l}^{2}{x_{{3
}}}^{2}{x_{{12}}}^{2}-4\,klx_{{3}}{x_{{12}}}^{2}-2\,{l}^{2}x_{{3}}{x_{
{12}}}^{2}-l{x_{{3}}}^{2}{x_{{12}}}^{2}+6\,{k}^{2}{x_{{3}}}^{2}\\&+4\,{k}
^{2}x_{{3}}x_{{12}}+3\,kl{x_{{3}}}^{2}+{l}^{2}{x_{{12}}}^{2}+4\,lx_{{3
}}{x_{{12}}}^{2}-12\,{k}^{2}x_{{3}}-6\,klx_{{3}}-7\,k{x_{{3}}}^{2}-5\,
kx_{{3}}x_{{12}}-2\,l{x_{{3}}}^{2}\\&-2\,l{x_{{12}}}^{2}+3\,kl+20\,kx_{{3
}}+4\,x_{{3}}l+2\,{x_{{3}}}^{2}+2\,x_{{12}}x_{{3}}-6\,k-2\,l-8\,x_{{3}
}+4=0
\end{aligned}\right.\]
Consider the ideal $I$ generated by $\{g_1, g_2, z x_3 x_{12}-1\}$ and the polynomial ring $R$ with coefficients in $\mathbb{Q}$, and take a lexicographic order $>$ with $z > x_3 > x_{12}$ for a monomial ordering on $R$. By the computer software, we get two polynomials containing in the Gr\"{o}bner basis for the ideal $I$.\\
$h(x_{12})={l}^{2} ( k+l )  ( 2\,{k}^{2}+2\,kl+{l}^{2}-l )  {x_{{12}}}^{8}\\-2\,{l}^{2} ( 2\,k+l-2 )  ( 4\,{k}^{2}+4\,kl+{l}^{2}-l
 )  {x_{{12}}}^{7}\\+ l ( 16\,{k}^{4}+76\,{k}^{3}l+71\,{k}^{2}{l}^{2}+22\,k{l}^{3}+{l}^
{4}-20\,{k}^{3}-119\,{k}^{2}l-81\,k{l}^{2}-16\,{l}^{3}+8\,{k}^{2}+51\,
kl+19\,{l}^{2}-4\,l )  {x_{{12}}}^{6}\\-4\,l ( 2\,k+l-2 )  ( 14\,{k}^{3}+20\,{k}^{2}l+5\,k{l}
^{2}-14\,{k}^{2}-21\,kl-4\,{l}^{2}+4\,k+4\,l ) {x_{{12}}}^{5}\\+ ( 32\,{k}^{5}+
344\,{k}^{4}l+368\,{k}^{3}{l}^{2}+117\,{k}^{2}{l}^{3}+6\,k{l}^{4}-80\,
{k}^{4}-842\,{k}^{3}l-686\,{k}^{2}{l}^{2}-168\,k{l}^{3}-4\,{l}^{4}+82
\,{k}^{3}+713\,{k}^{2}l+406\,k{l}^{2}+60\,{l}^{3}-40\,{k}^{2}-236\,kl-
76\,{l}^{2}+8\,k+20\,l ) {x_{{12}}}^{4}\\-2\, ( 2\,k+l-2 )  ( 48\,{k}^{4}+124\,{k}^{3}l+31\,{k}
^{2}{l}^{2}-92\,{k}^{3}-215\,{k}^{2}l-46\,k{l}^{2}+64\,{k}^{2}+110\,kl
+16\,{l}^{2}-16\,k-16\,l )  {x_{{12}}}^{3}\\+ ( 448\,{k}^{5}+608\,{k}^{4}l+212\,{k}
^{3}{l}^{2}+9\,{k}^{2}{l}^{3}-1424\,{k}^{4}-1550\,{k}^{3}l-448\,{k}^{2
}{l}^{2}-12\,k{l}^{3}+1714\,{k}^{3}+1427\,{k}^{2}l+304\,k{l}^{2}+4\,{l
}^{3}-956\,{k}^{2}-556\,kl-64\,{l}^{2}+232\,k+76\,l-16 ) {x_{{12
}}}^{2}\\-4\, ( k-1 )  ( 5\,k-2 )  ( -2+3\,k
 )  ( 2\,k+l-2 )  ( 4\,k+l-1 ) x_{{12}}\\+4\, ( 5\,k-2 ) ^{2} ( k-1 ) ^{2} ( l-1+2\,k) $\\
and\\
$h(x_{12},x_3)= -{l}^{2} ( k+l )  ( {k}^{2}+4\,kl+2\,{l}^{2}-2\,k-2\,l
 )  ( 2\,{k}^{2}+2\,kl+{l}^{2}-l ) {x_{{12}}}^{7}\\+ 2\,{l}^{2} ( 2\,k+l-2 )  ( {k}^{2}+4\,kl+2\,{l}^{2}-2
\,k-2\,l )  ( 4\,{k}^{2}+4\,kl+{l}^{2}-l ) {x_{{12}}}^{6}\\-l ( 16\,{k}^{6}+120\,{k}^{5}l+357\,{k}^{4}{l}^{2}+408\,{k}^{3}{l
}^{3}+206\,{k}^{2}{l}^{4}+43\,k{l}^{5}+2\,{l}^{6}-52\,{k}^{5}-365\,{k}
^{4}l-841\,{k}^{3}{l}^{2}-714\,{k}^{2}{l}^{3}-252\,k{l}^{4}-32\,{l}^{5
}+48\,{k}^{4}+357\,{k}^{3}l+622\,{k}^{2}{l}^{2}+355\,k{l}^{3}+66\,{l}^
{4}-16\,{k}^{3}-122\,{k}^{2}l-154\,k{l}^{2}-44\,{l}^{3}+8\,kl+8\,{l}^{
2} ) {x_{{12}}}^{5}\\+2\,l ( 2\,k+l-2 )  ( 28\,{k}^{5}+112\,{k}^{4}l+166\,{k
}^{3}{l}^{2}+90\,{k}^{2}{l}^{3}+15\,k{l}^{4}-84\,{k}^{4}-254\,{k}^{3}l
-278\,{k}^{2}{l}^{2}-114\,k{l}^{3}-14\,{l}^{4}+64\,{k}^{3}+172\,{k}^{2
}l+131\,k{l}^{2}+28\,{l}^{3}-16\,{k}^{2}-32\,kl-14\,{l}^{2} )  {x_{{12}}}^{4}\\+
 ( -32\,{k}^{7}-372\,{k}^{6}l-1118\,{k}^{5}{l}^{2}-1337\,{k}^{4}{
l}^{3}-710\,{k}^{3}{l}^{4}-153\,{k}^{2}{l}^{5}-7\,k{l}^{6}+144\,{k}^{6
}+1664\,{k}^{5}l+3882\,{k}^{4}{l}^{2}+3620\,{k}^{3}{l}^{3}+1499\,{k}^{
2}{l}^{4}+247\,k{l}^{5}+6\,{l}^{6}-242\,{k}^{5}-2641\,{k}^{4}l-4816\,{
k}^{3}{l}^{2}-3415\,{k}^{2}{l}^{3}-1003\,k{l}^{4}-98\,{l}^{5}+204\,{k}
^{4}+1882\,{k}^{3}l+2625\,{k}^{2}{l}^{2}+1293\,k{l}^{3}+210\,{l}^{4}-
88\,{k}^{3}-588\,{k}^{2}l-586\,k{l}^{2}-150\,{l}^{3}+16\,{k}^{2}+56\,k
l+32\,{l}^{2} ) {x_{{12}}}^{3}\\+2\, ( 2\,k+l-2 )  ( 48\,{k}^{6}+156\,{k}^{5}l+263\,{k}
^{4}{l}^{2}+162\,{k}^{3}{l}^{3}+27\,{k}^{2}{l}^{4}-188\,{k}^{5}-583\,{
k}^{4}l-706\,{k}^{3}{l}^{2}-342\,{k}^{2}{l}^{3}-48\,k{l}^{4}+248\,{k}^
{4}+688\,{k}^{3}l+631\,{k}^{2}{l}^{2}+216\,k{l}^{3}+20\,{l}^{4}-144\,{
k}^{3}-316\,{k}^{2}l-216\,k{l}^{2}-40\,{l}^{3}+32\,{k}^{2}+48\,kl+20\,
{l}^{2} ) {x_{{12}}}^{2}\\+ ( -368\,
{k}^{7}-740\,{k}^{6}l-850\,{k}^{5}{l}^{2}-513\,{k}^{4}{l}^{3}-120\,{k}
^{3}{l}^{4}-3\,{k}^{2}{l}^{5}+1908\,{k}^{6}+3454\,{k}^{5}l+3330\,{k}^{
4}{l}^{2}+1642\,{k}^{3}{l}^{3}+321\,{k}^{2}{l}^{4}+8\,k{l}^{5}-3740\,{
k}^{5}-5875\,{k}^{4}l-4684\,{k}^{3}{l}^{2}-1787\,{k}^{2}{l}^{3}-256\,k
{l}^{4}-4\,{l}^{5}+3546\,{k}^{4}+4678\,{k}^{3}l+2943\,{k}^{2}{l}^{2}+
784\,k{l}^{3}+60\,{l}^{4}-1676\,{k}^{3}-1786\,{k}^{2}l-800\,k{l}^{2}-
116\,{l}^{3}+344\,{k}^{2}+280\,kl+68\,{l}^{2}-16\,k-8\,l ) x_{{
12}}\\+2\, ( k-1 )  ( k-2 )  ( 5\,k-2 )
 ( 3\,{k}^{2}+2\,kl+{l}^{2}-2\,k-l )  (8\,{k}^{2}+4\,kl-8\,k ) \\+2\, ( l-1 )  ( k-1 )  ( k-2 )
 ( 5\,k-2 )  ( 2\,k+l-1 )  ( 3\,{k}^{2}+(2\,k+l)(l-1))x_3.$

By $h(x_{12},x_3)=0$, $x_3$ can be expressed by a polynomial of $x_{12}$ of degree $7$ with coefficients in $\mathbb{Q}$ for $k\geq3$ and $l\geq2$. It means if $x_{12}=s\in \mathbb R$ is a solution of $h(x_{12})=0$, then there exists $x_3=t\in\mathbb R$ such that $h(s,t)=0.$ Moreover, we have
 \begin{equation*}
 \begin{split}
 &h(0)=4\, ( 5\,k-2 ) ^{2} ( k-1 ) ^{2} ( l-1+2\,k
 ),\\
 &h(1)= ( k-1 )  ( l-1+2\,k )  ( k-l )
 ( 2\,k+l ) ^{2}\\
&h(+\infty)\rightarrow+\infty.
 \end{split}
 \end{equation*}
 For $l>k\geq 2$, we have $h(0)>0$ and $h(1)<0$. as a result, $h(x_{12})=0$ has at least two solutions, one of which is between $0$ and $1$ and the other is more than $1$. As is shown above, there exists $x_3\in\mathbb R$ with respect to each solution of $h(x_{12})=0$.

 The following is to check when $x_3\in \mathbb R^+$. For this, take a lexicographic order $>$ with $z > x_{12} > x_{3}$ for a monomial ordering on $R$. Similarly, we have the following polynomial containing in the Gr\"{o}bner basis for the ideal $I$:\\
 $p(x_3)=4\, ( l-1+2\,k )  ( 2\,{k}^{2}+2\,kl+l ( l-1
 )  )  ( 3\,{k}^{2}+2\,kl+{l}^{2}-2\,k-l ) ^{2}x_3^8\\
 -16\, ( 2\,k+l-2 )  ( 3\,{k}^{2}+ ( 2\,k+l
 )  ( l-1 )  ) (3\,{k}^{4}+ ( 12\,l-2 ) {k}^{3}+ ( 14\,{l}^{2}-10\,l
 ) {k}^{2}+ ( {l}^{2} ( 8\,l-10 ) +2\,l )
k+2\,{l}^{3} ( l-2 ) +2\,{l}^{2})x_3^7\\+(160\,{k}^{7}+ ( 1928\,l-1000 ) {k}^{6}+ ( 5492\,{l}^{2
}-6764\,l+1408 ) {k}^{5}+ ( 7644\,{l}^{3}-15830\,{l}^{2}+
7818\,l-736 ) {k}^{4}+ ( 6201\,{l}^{4}-18266\,{l}^{3}+15881
\,{l}^{2}-3848\,l+128 ) {k}^{3}+ ( 3072\,{l}^{5}-11752\,{l}
^{4}+14848\,{l}^{3}-6912\,{l}^{2}+808\,l ) {k}^{2}+ ( 880\,
{l}^{6}-4208\,{l}^{5}+7012\,{l}^{4}-4824\,{l}^{3}+1204\,{l}^{2}-64\,l
 ) k+112\,{l}^{7}-664\,{l}^{6}+1432\,{l}^{5}-1352\,{l}^{4}+504\,
{l}^{3}-32\,{l}^{2})x_3^6\\-4(2k+l-2)( ( 128\,l-108 ) {k}^{5}+ ( 582\,{l}^{2}-725\,l+144
 ) {k}^{4}+ ( 921\,{l}^{3}-1857\,{l}^{2}+800\,l-48 )
{k}^{3}+ ( 736\,{l}^{4}-2072\,{l}^{3}+1600\,{l}^{2}-324\,l
 ) {k}^{2}+ ( 312\,{l}^{5}-1152\,{l}^{4}+1292\,{l}^{3}-484
\,{l}^{2}+48\,l ) k+56\,{l}^{6}-268\,{l}^{5}+424\,{l}^{4}-252\,{
l}^{3}+40\,{l}^{2})x_3^5\\+( ( 200\,l-320 ) {k}^{6}+ ( 2448\,{l}^{2}-4584\,l+1712
 ) {k}^{5}+ ( 6812\,{l}^{3}-18568\,{l}^{2}+13190\,l-2216
 ) {k}^{4}+ ( 8618\,{l}^{4}-32023\,{l}^{3}+36747\,{l}^{2}-
13648\,l+1056 ) {k}^{3}+ ( 5781\,{l}^{5}-27998\,{l}^{4}+
45785\,{l}^{3}-29060\,{l}^{2}+6104\,l-160 ) {k}^{2}+ ( 2000
\,{l}^{6}-12260\,{l}^{5}+27136\,{l}^{4}-25936\,{l}^{3}+9932\,{l}^{2}-
1120\,l ) k+280\,{l}^{7}-2120\,{l}^{6}+6124\,{l}^{5}-8336\,{l}^{
4}+5284\,{l}^{3}-1280\,{l}^{2}+64\,l)x_3^4\\-(2(l-2))(2k+l-2)( ( 128\,l-108 ) {k}^{4}+ ( 582\,{l}^{2}-725\,l+144
 ) {k}^{3}+ ( 842\,{l}^{3}-1694\,{l}^{2}+716\,l-48 )
{k}^{2}+ ( 512\,{l}^{4}-1560\,{l}^{3}+1232\,{l}^{2}-244\,l
 ) k+112\,{l}^{5}-480\,{l}^{4}+632\,{l}^{3}-280\,{l}^{2}+32\,l)x_3^3\\+2(l-2)^2(20\,{k}^{5}+ ( 241\,l-125 ) {k}^{4}+ ( 659\,{l}^{2}-
800\,l+158 ) {k}^{3}+ ( 703\,{l}^{3}-1593\,{l}^{2}+820\,l-
68 ) {k}^{2}+ ( 328\,{l}^{4}-1156\,{l}^{3}+1168\,{l}^{2}-
336\,l+8 ) k+56\,{l}^{5}-276\,{l}^{4}+448\,{l}^{3}-268\,{l}^{2}+
48\,l)x_3^2\\-2\, ( l-2 ) ^{3} ( 2\,k+l-2 )  ( -2+3\,k+2
\,l )  ( 3\,{k}^{2}+2k(6l-1)+4l(2l-1) )x_3 \\+( l-2 ) ^{4} ( k+l )  ( -2+3\,k+2\,l) ^{2}.$

It is not hard to verify that the coefficients of polynomial $p(x_3)$ is positive for even degree and negative for odd degree whenever $l>k\geq3$. That is, all the solution of $x_3$ for homogeneous Einstein equations are positive (if exist) whenever $l>k\geq3$.

In summary, for $l>k\geq3$, we get two solutions of the form
$$\{x_1=x_2={\frac {x_{{12}} \left( k-2 \right) }{l{x_{{12}}}^{2}+3\,k-2}}, x_3=\alpha(x_{12}), x_{12}\neq1, x_{13}=x_{23}=1\},$$
where $\alpha(x_{12})$ is a rational polynomial of $x_{12}$ with positive values. By Proposition \ref{prop}, every Einstein metric induced by the solution of this form is non-naturally reductive. Thus, we have proved Theorem \ref{thm2}.

By Theorem~\ref{thm2}, for any $3\leq k\leq [\frac{n-1}{3}]$, $\SO(n)$ admits at least two non-naturally reductive left-invariant Einstein metrics which are $\Ad(\SO(k)\times\SO(k)\times\SO(n-2k))$-invariant. That is, Theorem~\ref{thm} follows.


\begin{thebibliography}{10}
\bibitem{ArMoSa}
A. Arvanitoyeorgos, K. Mori and Y. Sakane, \emph{Einstein metrics on compact Lie groups which are not naturally reductive}, Geom. Dedicate 160(1) (2012), 261--285.

\bibitem{ArSaSt3}
A. Arvanitoyeorgos, K. Mori and Y. Sakane, \emph{New Einstein metrics on the Lie group $SO(n)$ which are not naturally reductive}, arXiv: 1511.08849v1, 2015.

\bibitem{Be}
A.L. Besse, \emph{Einstein Manifolds}, Springer-Verlag, Berlin, 1986.

\bibitem{Bo}
C. Bohm, \emph{Homogeneous Einstein metrics and simplicial complexes}, J. Differential Geom. 67(1) (2004), 74--165.

\bibitem{BoWaZi}
C. Bohm, M. Wang and W. Ziller, \emph{A variational approach for compact homogeneous Einstein manifolds}, Geom. Func. Anal. 14(4) (2004), 681--733.

\bibitem{CC}
H. Chen and Z. Chen, \emph{Notes on ``Einstein metrics on compact simple Lie groups attached to standard triples"}, arXiv: 1701.01713v1, 2017.

\bibitem{ChLi}
Z. Chen and K. Liang, \emph{Non-naturally reductive Einstein metrics on the compact simple Lie group $F_4$}, Ann. Glob. Anal. Geom. 46 (2014), 103--115.

\bibitem{ChSa}
I. Chrysikos and Y. Sakane, \emph{Non-naturally reductive Einstein metrics on exceptional Lie groups}, arXiv: 1511.03993v1, 2015.

\bibitem{JW} J.E. D' Atri and W. Ziller, \emph{Naturally reductive metrics and Einstein metrics on compact Lie groups}, Memoirs Amer. Math. Soc. 19 (215), 1979.

\bibitem{He}
J. Heber, \emph{Noncompact homogeneous Einstein spaces}, Invent. math. 133 (1998), 279--352.

\bibitem{La}
J. Lauret, \emph{Einstein solvmanifolds are standard}, Ann. Math. (2) 172(3) (2010), 1859--1877.

\bibitem{Mo}
K. Mori, \emph{Left Invariant Einstein Metrics on $SU(n)$ that are not naturally reductive}, Master Thesis (in Japanese) Osaka University 1994, English Translation: Osaka University RPM 96010 (preprint series) 1996.

\bibitem{PaSa}
J.S. Park and Y. Sakane, \emph{Invariant Einstein metrics on certain homogeneous spaces}, Tokyo J. Math. 20(1) (1997), 51--61.

\bibitem{W1} M. Wang, \emph{Einstein metrics from symmetry and bundle constructions}, In:  Surveys in Differential Geometry: Essays on Einstein Manifolds. Surv. Differ. Geom. VI, Int. Press, Boston, Ma 1999.

\bibitem{W2} M. Wang, \emph{Einstein metrics from symmetry and bundle constructions: A sequel}, In:  Differential Geometry: Under the Influence of S.-S. Chern, in: Advanced Lectures in Mathematics, vol. 22, Higher Education Press/International Press, 2012, pp. 253--309.

\bibitem{WaZi} M. Wang and W. Ziller, \emph{Existence and non-existence of homogeneous Einstein metrics}, Invent. Math. 84 (1986), 177--194.

\bibitem{YD1}
Z. Yan and S. Deng, \emph{Einstien metrics on compacr simple Lie groups attached to standard triples}, to appear in Trans. Amer. Math. Soc., 2016.





\end{thebibliography}
\end{document}